\documentclass[10pt,conference,letterpaper]{IEEEtran}

\IEEEoverridecommandlockouts

\usepackage{times,amsmath,epsfig,amssymb,graphicx,amsfonts,amsthm}
\usepackage{subfigure}
\usepackage{empheq}
\usepackage{psfrag}
\usepackage{fge}
\usepackage{bbm}
\usepackage{tikz}
\usepackage{mathtools}
\usepackage{mathdots}
\usepackage{amsmath}
\usepackage{relsize}
\mathtoolsset{showonlyrefs}
\usepackage{amsmath}
\usepackage{amsthm}
\usepackage{amsfonts}   
\usepackage{amssymb}
\usepackage{mathrsfs}
\usepackage{setspace}
\usepackage{algorithm}
\usepackage{algcompatible}
\usepackage{enumerate}
\usepackage{stmaryrd}
\usepackage{color}
\usepackage{xcolor}
\usepackage{soul}

\floatname{algorithm}{Routine}

\newtheorem{mydef}{Definition}
\newtheorem{mytheorem}{Theorem}
\newtheorem{myassumption}{Assumption}
\newtheorem{mylemma}{Lemma}

\newtheorem{mycorollary}{Corollary}
\newtheorem{myproposition}{Proposition}

\newtheorem{myproblem}{Problem}

\newcounter{ale}

\newenvironment{liste}{\begin{itemize}}{\end{itemize}}
\newcommand{\aliste}{\begin{liste} \setcounter{ale}{1}}
\newcommand{\zliste}{\end{liste}}

\newcommand{\bigslant}[2]{{\raisebox{.2em}{$#1$}\left/\raisebox{-.2em}{$#2$}\right.}}

\title{{{\bf  \LARGE  Failure Detection and Isolation in Integrator Networks}}}

\author{Mohammad Amin Rahimian, Victor M. Preciado{\small $~^{*}$}
\thanks{$^{*}$ The authors are with the Department of Electrical and Systems Engineering, University of Pennsylvania, Philadelphia, PA 19104-6228 USA. (email: {\fontsize{8}{8}\selectfont\ttfamily\upshape preciado@seas.upenn.edu}).}
}

\begin{document}

\maketitle
\begin{abstract}
Detection and isolation of link failures under the Laplacian consensus dynamics have been the focus of our previous study. Our results relate the failure of links in the network to jump discontinuities in the derivatives of the output responses of the nodes and exploit that relation to propose failure detection and isolation (FDI) techniques, accordingly. In this work, we extend the results to general linear networked dynamics. In particular, we show that with additional niceties of the integrator networks and the enhanced proofs, we are able to incorporate both unidirectional and bidirectional link failures. At the next step, we extend the available FDI techniques to accommodate the cases of bidirectional link failures and undirected topologies. Computer experiments with large networks and both directed and undirected topologies provide interesting insights as to the role of directionality, as well as the scalability of the proposed FDI techniques with the network size.
\end{abstract}

\section{Introduction}
Multi-agent network systems have found promising applications in areas such as motion coordination of robots \cite{mesbahiBook}. Cooperative dynamics over a network can be strongly affected by the network failures. Hence, studying the effects of link or node failures on the network dynamics is an important topic in network science and it has various practical implications \cite{aminAutomatica,Kleinberg}. 

Fault Detection and Isolation (FDI) in networked systems is an active area of research with application to power networks \cite{largeScaleNonlinearPowerNetworks}, and security of cyber-physical systems \cite{6545301}. In a systematic approach the so-called FDI filter uses available measurements to generate a residual signal, which is then used for fault diagnosis \cite{discreteTimeFaultDetection}. Various observer and Kalman filter techniques are used to obtain the residual signal \cite{6567892}, which can be then compared with a threshold signal to detect faults. In \cite{6482175} the residual signal and threshold signals are generated using filtering techniques that allow for noise suppression and less conservative detection thresholds. Failures and attacks are modeled as disturbances in the descriptor systems approach that is adopted in \cite{6545301}, and some fundamental limitations of detection and identification in cyber-physical systems are established from systems-theoretic or graph-theoretic perspectives.

The mathematical framework for investigating the ability to detect and distinguish failures based on the observed output responses is established in \cite{asjc}, where conditions in terms of the inter-nodal distances to the observation points are provided for the detectability of links.  Subsequent results in \cite{CDC13} provide a method for detection and isolation of failures under the Laplacian network dynamics. These results relate the presence of discontinuities in the derivatives of particular orders in the output responses of a subset of nodes to occurrence and location of link failures; and coupled with efficient sensor placement algorithms they can be used to pin down the location of failures in the network. In this paper, we analyze the particular case of networked single integrator dynamics and extend the proofs and FDI algorithms to allow for efficient sensor location in a (directed or undirected) integrator network where link failures can be either unidirectional or bidirectional.

The remainder of this paper is organized as follows. In Section~\ref{sec:pre} we set up the notation and give the preliminaries on weighted digraphs and their algebraic properties. Our main analytical findings are set forth in Section~\ref{sec:singleIntegrator} where we first prove the results relating link failures to jump discontinuities of the output derivatives for single-integrator networks, and then  consider the case when the links can fail simultaneously in both forward and reverse directions. In Section~\ref{sec:singleIntegratorSensorPlacement} we flesh out the necessary formalism to accommodate both cases of unidirectional and bidirectional link failures within the framework of the previously developed sensor placement algorithms. Computer experiments on large networks in Section~\ref{sec:examples} elucidate the results and offer interesting insights on the effects of the network size and edge directionality. Section~\ref{sec:conc} concludes the paper.

\section{Preliminaries}\label{sec:pre}

Throughout the paper, $\varnothing$ is the empty set, $\mathbb{N}$ denotes the set of all natural numbers, and $\mathbb{R}$ denotes the set of all real numbers. Also, the set of integers $\{1,2,\ldots,k\}$ is denoted by $\mathbb{N}_k$, and any other set is represented by a calligraphic capital letter. The cardinality of a set $\mathcal{X}$ is denoted by $|\mathcal{X}|$, and $\mathcal{P}(\mathcal{X}) = \{ \mathcal{M}; \mathcal{M} \subset \mathcal{X} \}$ denotes the power-set of $\mathcal{X}$, which is the set of all its subsets. The difference of two sets $\mathcal{X}$ and $\mathcal{Y}$ is  denoted by $\mathcal{X} \fgebackslash \mathcal{Y}$ and is defined as $\left\{x;x \in \mathcal{X} \wedge x \notin \mathcal{Y}\right\}$, where $\wedge$ is the logical conjunction. In addition the logical implication and bi-implication are denoted by $\rightarrow$ and $\leftrightarrow$, respectively, and $\vee$ denote the logical disjunction. The fixed integer $n$ represents the number of nodes in the network. Matrices are represented by capital letters, vectors are expressed by boldface lower-case letters, and the superscript $^{T}$ denotes the matrix transpose. Moreover, $I$ denotes the identity matrix with proper dimension, $\mathbf{e}_i$ is the $i$-th unit vector in the standard basis of $\mathbb{R}^n$, and $m_{ij}:=\left[M\right]_{ij}$ indicate the element of matrix $M$ which is located at its $i-$th row and $j-$th column. 

A directed graph or \emph{digraph} is defined as an ordered pair of sets $\mathcal{G} := (\mathcal{V},\mathcal{E})$, where $\mathcal{V} = \left\{{\nu}_1,\ldots,{\nu}_n\right\}$ is a set of $n = |\mathcal{V}|$ vertices and $\mathcal{E} \subseteq \mathcal{V} \times \mathcal{V}$ is a set of directed edges. In the graphical representations, each edge $\epsilon := ({\tau},{\nu}) \in \mathcal{E}$ is depicted by a directed arc from vertex ${\tau} \in \mathcal{V}$ to vertex ${\nu} \in \mathcal{V}$. Vertices ${\nu}$ and $\tau$ are referred to as the \emph{head} and \emph{tail} of the edge $\epsilon$ and a $(\nu,\nu)$ edge is called a self-loop on $\nu$.   Given an integer $k \in \mathbb{N}$, a set of (possibly repeated) indices $\{\alpha_1,\alpha_2,\ldots,\alpha_k\}$ $ \subseteq$ $\mathbb{N}_{n}$ and two vertices $\tau,\nu $ $\in$ $ \mathcal{V}$, an ordered sequence of edges of the form $\mathcal{W}$ $:=$ $({\tau},{\nu}_{\alpha_1}),({\nu}_{\alpha_1},{\nu}_{\alpha_2}),$ $ \ldots$ $,$ $({\nu}_{\alpha_{k-1}},{\nu}_{\alpha_{k}}),({\nu}_{\alpha_{k}},{\nu})$ is called a \emph{${\tau}{\nu}$ walk} with start-node $\tau$, end-node $\nu$ and \emph{length} $k+1$. A cycle on node $\nu$ signifies a $\nu\nu$ walk. For any $\{{q},{p}\} \subset \mathbb{N}_{n}$, $\Omega^{k}(\mathcal{G};{\nu_q,\nu_p})$ is the set of all ${\nu_q}{\nu_p}$ walks in $\mathcal{G}$ with length $k$. Similarly for $\{{q},{r},{p}\}$ $\subset$ $\mathbb{N}_{n},$ $\Omega^{k}(\mathcal{G};{\nu_q,\nu_r,\nu_p})$ $=$ $\{ \mathcal{W}$ $\in$ $\Omega^{k}(\mathcal{G};{\nu_q,\nu_p}); (\nu_q,\nu_r)$ $\in$ $\mathcal{W} \}$, i.e. the set of $\nu_q \nu_p$ walks with  length $k$ that include the edge $(\nu_q,\nu_r)$. In the same venue, the integer $\textrm{dist}(\mathcal{G};{\nu_q},{\nu_p}) = \min\{k\in\mathbb{N}:,\Omega^{k}(\mathcal{G};{\nu_q,\nu_p}) \neq \varnothing\}$ is referred to as the distance from ${\nu_q}$ to ${\nu_p}$ in $\mathcal{G}$, and by convention, $\textrm{dist}(\mathcal{G};{\nu_q},{\nu_q}) = 0, \forall q \in \mathbb{N}_n$ and $\textrm{dist}(\mathcal{G};{\nu_q},{\nu_p}) = \infty$ if $\forall k\in\mathbb{N}$, $\Omega^{k}(\mathcal{G};{\nu_q,\nu_p}) = \varnothing$. For ease of notation we usually use $\textrm{dist}({\nu_q},{\nu_p}) := \textrm{dist}(\mathcal{G};{\nu_q},{\nu_p})$. The diameter of $\mathcal{G}$, denoted by $\mbox{diam}(\mathcal{G})$, is defined as the maximum distance between any pair of nodes $\nu_{p},\nu_{q}\in\mathcal{V}$: $\mbox{diam}(\mathcal{G}) = \max_{\nu_{q},\nu_{p} \in\mathcal{V}} \textrm{dist}(\nu_{q},\nu_{p})$.
A matrix $A \in \mathbb{R}^{n \times n}$ is called an in-weighting on $\mathcal{G}$ if $\forall \{ \nu_p,\nu_q\} \subset \mathcal{V}, (\nu_q,\nu_p) \not\in \mathcal{E} \rightarrow a_{pq} = 0$. If $A$ is an in-weighting on $\mathcal{G}$, then $\omega(\mathcal{W},A) = \prod _{(\nu_q,\nu_p)\in\mathcal{W}} a_{pq}$ is referred to as the weight of walk $\mathcal{W}$ w.r.t. $A$. Given a set of walks $\Omega$ in digraph $\mathcal{G}$ and the in-weighting $A$ on $\mathcal{G}$, the function $\Phi(\Omega, A) = \sum_{\mathcal{W}\in\Omega} \omega(\mathcal{W},A)$ is defined, which finds use in the proof of the main results in Section~\ref{sec:singleIntegrator}. Note that this function satisfies $\Phi(\Omega^{k}(\mathcal{G};{\nu_q,\nu_r,\nu_p}), A) = a_{rq} \Phi(\Omega^{k-1}(\mathcal{G};{\nu_r,\nu_p}), A)$. It is also known that given an in-weighting $A$ on $\mathcal{G}$ and vertices $\{\nu_q,\nu_p\} \subset \mathcal{V}$, that \cite{BiggsGraphTheory,Preciado:2013:MSA:2502376.2502379}:
\begin{align}
&\Phi(\Omega^{k}(\mathcal{G};{\nu_q,\nu_p}),A) = \left[A^{k}\right]_{pq}.
\label{eq:numWalks}
\end{align} 

The following definition comes handy in the statement and proof of the main theorem in the next section. At each time $t_f$ and for all $p \in \mathbb{N}_n$ and $ k \in \mathbb{N}_z$, define $\Delta_{t_f}: \mathbb{N}_n \times \mathbb{N}_z \to \mathbb{R}$ as $\Delta_{t_f}(p,k) := \frac{ {\textrm{d}}^{k} } {\textrm{d} t^{k}}({x}_p)(t_f^{+}) - \frac{{\textrm{d}}^{k} }{\textrm{d} t^{k}}(x_p)(t_f^{-})$. The function $\Delta_{t_f}(p,k)$ as defined, measures the jump in the $k-$th derivative of the response of agent $p$ at the time of failure $t_f$, and the parameter $z\leqslant \mbox{diam}(\mathcal{G})+1$ is a fixed integer, denoting the highest order of derivatives to which the designer has access. 

\section{Failures in Networks of Single-Integrators}\label{sec:singleIntegrator}

Let us consider a network of $n$ single-integrators, where each integrator is described by a single state $x_{i}$, with the following dynamics:
\begin{equation}
\dot{\mathbf{x}}(t)=A\mathbf{x}(t)+B\mathbf{w}(t),\: t>t_{0},\label{eq:laplaciandynamics}
\end{equation}
where $\mathbf{x}\left(t\right)=\left(x_{1},\ldots,x_{n}\right)^{T}\in\mathbb{R}^{n}$,
$\mathbf{w}(t)\in\mathbb{R}^{m}$, $B\in\mathbb{R}^{N\times m}$, and $A\in\mathbb{R}^{N\times N}$ is an in-weighting of graph $\mathcal{G}$. We assume that the entries of $\mathbf{w}\left(t\right)$ are $\left(\mathrm{diam}(\mathcal{G})+1\right)$-differentiable. Let us assume link $\bar{\epsilon}=(\nu_{j},\nu_{i})$ fails at time
$t=t_{f}$, resulting in a faulty connectivity graph $\bar{\mathcal{G}}=(\mathcal{V},\mathcal{E}\fgebackslash\{\bar{\epsilon}\})$
for $t>t_{f}$. The corresponding in-weighting of $\bar{\mathcal{G}}$,
denoted by $\bar{A}$, is a perturbed version of $A$ that satisfies $\bar{a}_{ij}=0$, and $\bar{a}_{qr}=a_{qr}$ for all $q\neq i$; i.e. entries on the $i$-th row of $A$ are allowed to change while every other entry remains unaffected.

The following theorem characterizes the effect of link failures on the output derivatives of a network of single integrators:

\begin{mytheorem}\label{theo:detection2} Consider the dynamic
network of single integrators in (\ref{eq:laplaciandynamics}) with
output equation $y_{p}(t)=x_{p}(t)$, and assume link $\bar{\epsilon}=(\nu_{j},\nu_{i})$
fails at time $t_{f}$. Then, the output of node $\nu_{p}$ satisfies:
\begin{align}
\Delta_{t_f}({p,k})=\Phi(\Omega^{k-1}(\nu_{i},\nu_{p}),A)\sum_{q=1}^{N}\left(\bar{a}_{iq}-a_{iq}\right)x_{q}(t_{f}),\label{eq:mainResult}
\end{align}
for $k=\mbox{\emph{dist}}(\nu_{j},\nu_{p})$; and $\Delta_{t_f}({p,k})=0$,
for $k<\mbox{\emph{dist}}(\nu_{j},\nu_{p})$. \end{mytheorem}

\begin{IEEEproof} See the Appendix. \end{IEEEproof}

The next lemma and the following corollary show in particular how the result of Theorem~\ref{theo:detection2} agrees with and maps to the result of Theorem~1 in \cite{rahimian2014detection}, where we consider the case of networked LTI agents. However, the proofs in \cite{rahimian2014detection} rely on Laplace domain techniques rather than the combinatorial arguments in the time-domains that helped us prove Theorem~\ref{theo:detection2}; indeed, the flexibilities in the latter case facilitate some extensions that are of particular interest to FDI scenarios.

\begin{mylemma}\label{lem:singleIntegratorDetection} Consider the failure of link $\bar{\epsilon} = (\nu_j,\nu_i)$ and suppose that $k=\mbox{\emph{dist}}(\nu_{j},\nu_{p})$. If the failure of link $\bar{\epsilon}$ is to induce a jump discontinuity in the $k$-th derivative of the output response of node $\nu_{p}$, then it should be true that $\mbox{\emph{dist}}(\nu_{i},\nu_{p})$ $=$ $k-1$ $<$ $\mbox{\emph{dist}}(\nu_{j},\nu_{p})$ $=$ $k$.
\end{mylemma}

\begin{IEEEproof} See the Appendix. \end{IEEEproof}
The following is now immediate upon combining the results of Theorem~\ref{theo:detection2} and  Lemma~\ref{lem:singleIntegratorDetection}.
\begin{mycorollary}\label{coro:detection4} Under the assumptions of Theorem~\ref{theo:detection2},
\begin{align}
\Delta_{t_f}({p,k})=\left\{ \begin{array}{ll}
{c}_{i,j,p}, & \mbox{for }k=(\mbox{\emph{dist}}(\nu_{i},\nu_{p})+1),\\
0, & \mbox{for }k<(\mbox{\emph{dist}}(\nu_{i},\nu_{p})+1),
\end{array}\right.\label{eq:limit3}
\end{align} where ${c}_{i,j,p} = \Phi(\Omega^{k-1}(\nu_{i},\nu_{p}),A)\sum_{q=1}^{N}\left(\bar{a}_{iq}-a_{iq}\right)x_{q}(t_{f})$.
\end{mycorollary}

The condition stated in Lemma~\ref{lem:singleIntegratorDetection} is intuitive because $\mbox{{dist}}(\nu_{i},\nu_{p})$ $=$ $k-1$ $<$ $\mbox{{dist}}(\nu_{j},\nu_{p})$ $=$ $k$ holds true only if there exist a shortest path of length $k$ connecting $\nu_{j}$ to $\nu_{p}$ and with $(\nu_{j},\nu_{i})$ as its first edge. In other words, the failed link $(\nu_{j},\nu_{i})$ contributes to the flow of information from $\nu_{j}$ to $\nu_{i}$ as an element of a shortest path from node $\nu_{j}$ to node $\nu_{i}$. These observations are in perfect agreement with the sufficient conditions previously studied in \cite{asjc}.  On the other hand, Corollary~\ref{coro:detection4} shows how Theorem~\ref{theo:detection2} may follow as a special case of Theorem~1 in \cite{rahimian2014detection} (up to a known constant multiplier), after setting the relative degree $r=1$ for the involved networked LTI systems. The proofs in the single-integrator case, however, admit additional niceties that we discuss next. 

\subsection{Bidirectional Link Failures}\label{sec:bidirectional}

To begin, note that the perturbed matrix $\bar{A}$ is not constrained in the way its entries on the $i$-th row are modified, thence Theorem~\ref{theo:detection2} continues to hold in the case where all or several of the edges incoming to node $\nu_i$ are lost simultaneously. Indeed, it is perceivable for the faults in an agent's hardware or internal structure to cause the failure of multiple links which are incoming to that agent. In the particular case of a faulty agent $\nu_i$, which looses all its incoming links at the instant of failure $t = t_f$, the systems dynamics for $t > t_f$ is characterized by $\bar{a}_{iq} = 0$, $\forall q \in \mathbb{N}_n\fgebackslash\{i\}$. Hence, as a special case, Theorem~\ref{theo:detection2} and the corresponding FDI techniques that are developed in Section~\ref{sec:singleIntegratorSensorPlacement} may also be applied for the detection and isolation of single agent failures by mapping the isolated edges to their head vertices.

On the other hand, for certain applications, where communications are of a bidirectional nature, it is reasonable to consider link failures that simultaneously prevent either agents from communication in the other's direction. Such a failure corresponds to the simultaneous elimination of both links $\bar{\epsilon}$, defined earlier, and $\hat{\epsilon} = (\nu_i,\nu_j)$ leading to $\hat{\mathcal{G}} = (\mathcal{V},\mathcal{E}\fgebackslash\{\bar{\epsilon},\hat{\epsilon}\})$ as the information flow structure for $t>t_f$. It is worth highlighting that undirected networks, where $\forall \{\tau,\nu\} \subset \mathcal{V}$, $(\tau,\nu) \in \mathcal{E} \leftrightarrow (\nu,\tau) \in \mathcal{E}$, signify the special case that all, not just some, of the links are bidirectional. The FDI methods in this paper are designed to handle the cases where some of the links in the networks are bidirectional and the rest are unidirectional. It is worth pinpointing that, as an assumption of modeling each link is considered either bidirectional or unidirectional, but not both. In other words, if $(\tau,\nu) \in \mathcal{E}$ and $(\nu,\tau) \in \mathcal{E}$ and the link between the nodes $\tau$ and $\nu$ is specified as bidirectional, then $(\tau,\nu)$ fails if, and only if, $(\nu,\tau)$ fails; otherwise, the two links $(\tau,\nu)$ and $(\nu,\tau)$ are regarded as separate, and their failures as independent events. Accordingly, it is assumed that the set of all bidirectional links $\mathcal{B}$ in the network, is known to the designer beforehand.

After modeling the failure of a bidirectional link as the simultaneous failures of two directed links, $\bar{\epsilon}=(\nu_{j},\nu_{i})$ and $\hat{\epsilon}=(\nu_{i},\nu_{j})$, the proof of Theorem~\ref{theo:detection2} can be adapted to yield:

\begin{myproposition}\label{prop:detection3} In the case of the simultaneous failure of $\bar{\epsilon}=(\nu_{j},\nu_{i})$ and $\hat{\epsilon}=(\nu_{i},\nu_{j})$,
\eqref{eq:limit3} still holds true if we substitute $\mbox{\emph{dist}}(\nu_{i},\nu_{p}) + 1$ by $\max\{\mbox{\emph{dist}}(\nu_{j},\nu_{p})+1$, $\mbox{\emph{dist}}(\nu_{i},\nu_{p})+1\}$. \end{myproposition}

\begin{IEEEproof} See the Appendix. \end{IEEEproof}

We end this section by an intuitive remark that as each agent of the network system in \eqref{eq:laplaciandynamics} is a single-integrator, a jump discontinuity (because of a sudden network failure) at point $\nu_i$ will appear to point $\nu_p$ after several (length of path) steps of ``integrations''. Thus, an agent at point $\nu_p$ needs to make the same number of ``differentiations'' before observing the jump due to the failure at point $\nu_i$. In what follows we shall see how to determine the observation points along with the required number of differentiations at each point so that the occurrence and location of failures are always inferable from the observed jumped discontinuities.

\section{Sensor Placement for Unidirectional and Bidirectional Links}\label{sec:singleIntegratorSensorPlacement}

It is assumed that at each instant of time, the designer is given access to the response of a subset of agents, as well as the nominal network information flow digraph $\mathcal{G}$ (prior to the link failure) and the set of bidirectional $\mathcal{B}$.  Neither the location of the failure (nodes $\nu_i$ and $\nu_j$), nor the time of failure $t_f$ are known to the designer. In the case of detection, the designer is interested in determining the existence of any single link failure in the network at the instant of failure. For the isolation problem, however, the designer would like to determine ``instantaneously", not only the existence of a failure, but also its location. That is to determine which link, if any, has failed and exactly at the same instant as it fails. The significance of \emph{``instantaneous"} detection and isolation is better understood upon noting that if the time of failure is random and has a continuous sample space, then  \emph{``simultaneous''} failure of more than one link is a measure zero event; hence, justifying the focus of investigation in this paper, which is on the \emph{``single"} (possibly bidirectional) link failures. Before shifting attention to the sensor placement problem, two assumptions are set forth:

\begin{myassumption}\label{assum:assumption1}For all pairs of nodes
$\nu_{p},\nu_{q}\in\mathcal{V}$, the in-weighting $A$ on digraph
$\mathcal{G}$ satisfies $\Phi(\Omega^{\mbox{\emph{dist}}(\nu_{q},\nu_{p})}(\nu_{q},\nu_{p}),A)$ $\neq$ $0$,
i.e., the sum of the weights of all shortest paths between them is
nonzero.\end{myassumption}

\begin{myassumption}\label{assum:assumption1-1} Given the in-weightings
$A$ and $\bar{A}$ of the faultless and the faulty network, $\mathcal{G}$
and $\bar{\mathcal{G}}$, respectively, we have
that $\sum_{q=1}^{N}\left(\bar{a}_{iq}-a_{iq}\right)x_{q}(t_{f})\neq0$, where $\nu_i$ is the head (or tail too, if the link is bidirectional) of the failed link and $t_{f}$ denotes the instant of failure.\end{myassumption}

The first assumption above is a provision of consistency that is assumed with regard to the in-weighting matrix $A$. This assumptions is satisfied almost surely for any assignment of weights on the graph. In particular,
it holds true for the Laplacian consensus networks considered in \cite{CDC13}. The second assumption involves the values of the agents states at the time failure $t_{f}$. This condition also holds true, almost surely, for
any in-weightings $A$, its perturbed version $\bar{A}$, and a random time of failure $t_{f}>t_{0}$; since $\sum_{q=1}^{N}\left(\bar{a}_{iq}-a_{iq}\right)x_{q}(t_{f})=0$ specifies a low-dimensional hyperplane in the agents' state space that the agents almost surely avoid given a random time of failure.

To enable the designer to handle the desired FDI tasks, she is given access to the output response of a subset of nodes as well as their derivatives upto the $z$-th order. In this section, we offer efficient procedures for determining such a subset of nodes, given the network topology and parameter $z$; and in such a way that all link failures in the networks can be detected or isolated from the occurrence of jump discontinuities in the observed outputs and their derivatives. Furthermore, we would like to achieve this goal using as few observation points (sensors) as possible. From Corollary~\ref{coro:detection4} it follows that if the existence of a jump discontinuity in the $k-$th derivative of the output response of agent $p$ is to serve as the basis for a method to detect the failure of edge $\bar{\epsilon}$ at time $t_f$, then it should be true that $\mbox{dist}({\mathcal{G}};{\nu}_i,{\nu_p})$ $ = $ $k-1$ $<$ $\mbox{dist}({\mathcal{G}};{\nu}_j,{\nu_p})$ $=$ $k$. The latter happens only if there exist a shortest path of length $k$ connecting $\nu_j$ to $\nu_p$ and with $(\nu_j,\nu_i)$ as its first edge. In \cite{CDC13} we use this observation in the case of Laplacian network dynamics to define binary relations  $\mathcal{R}_{k}, k \in \mathbb{N}_z$ and $\mathcal{R}_{0}$ between the sets $\mathcal{V}$ and $\mathcal{E}$ such that for all $p\in \mathcal{V}$ and $\epsilon \in \mathcal{E}$ if $(p,\epsilon) \in \mathcal{R}_{k}$, then the failure of link $\epsilon$ produces a jump in the $k-th$ derivative of the response of node $p$ and if $(p,\epsilon) \in \mathcal{R}_{0}$ then the failure of edge $\epsilon$ does not produce a jump in any of the derivatives of the response of node $p$ upto the $z-$th order. We now go ahead and redefine the binary relations per Proposition~\ref{prop:detection3} to accommodate bidirectional link failures. Indeed, bidirectional links are treated specially as any of the two edges in reverse directions can provide us with the required relation for detection when a bidirectional failure occurs. 

\begin{mydef} \label{Def:BinaryRelationships-1} We define the binary
relations $\mathcal{{R}}_{0}$ and ${\mathcal{R}}_{k}$
for $k\in\mathbb{N}_{z}$, between $\mathcal{V}$ and $\mathcal{E}$,
as follows. For all $p\in\mathcal{V}$ and $\epsilon=(\nu_{q},\nu_{r})\in\mathcal{E}$,
we have that:
\begin{itemize}
\item If $\epsilon\not\in\mathcal{B}$, then $(p,\epsilon)\in\mathcal{R}_{k}$ if, and only if, $\textrm{\emph{dist}}(\nu_{q},\nu_{p})=k\mbox{\emph{ and} }\textrm{\emph{dist}}(\nu_{r},\nu_{p})=k-1$.
\item If $\epsilon\in\mathcal{B}$, then $(p,\epsilon)\in\mathcal{R}_{k}$
if, and only if, one of the following conditions is satisfied:
\begin{itemize}
\item $\textrm{\emph{dist}}(\nu_{q},\nu_{p})=k\mbox{ \emph{and} }\textrm{\emph{dist}}(\nu_{r},\nu_{p})=k-1$,
or
\item $\textrm{\emph{dist}}(\nu_{r},\nu_{p})=k\mbox{ \emph{and} }\textrm{\emph{dist}}(\nu_{q},\nu_{p})=k-1$.
\end{itemize}
\end{itemize}
\end{mydef}

The FDI problems can now be posed as follows.

\begin{myproblem}[Detection]\label{prob:detection} Given a digraph $\mathcal{G}=(\mathcal{V},\mathcal{E})$,
find a subset of nodes $\mathcal{M}_{D}\subseteq\mathcal{V}$ of minimum cardinality $\left|\mathcal{M}_{D}\right|$, such that for all $\epsilon\in\mathcal{E}$, there exists a node ${p}\in\mathcal{M}_{D}$ such that $({p},\epsilon)\not\in\mathcal{R}_{0}$.
\end{myproblem}

\begin{myproblem}[Isolation]\label{prob:isolation} Given a digraph $\mathcal{G}=(\mathcal{V},\mathcal{E})$, find a subset of vertices $\mathcal{M}_{I}\subseteq\mathcal{V}$ with the smallest cardinality $\left|\mathcal{M}_{I}\right|$, such that $f_{I}\left(\mathcal{M}_{I}\right)=0$.
\end{myproblem}

The idea for proposing efficient sensor placement algorithms that approximate the solutions of the above problems is by counting the number of edges that are not yet detectable or isolatable from the currently chosen nodes and add new nodes to the existing sets in a greedy manner: in each addition of a new node to the existing sensor set, we aim to decrease the number of edges that are not yet detectable or isolatable as much as possible. To this end, we define a correspondence $\mathcal{I}:$ $\mathcal{P}(\mathcal{V})\times \mathcal{E}$ $\to$ $\mathcal{P}({(\mathbb{N}_z\cup\{0\}) \times \mathcal{V}})$ which maps an order pair $(\mathcal{M},\epsilon)$, comprised of a sensor set $\mathcal{M}$ and an edge $\epsilon$, to the set of ordered pairs that specify the relations between edge $\epsilon$ and nodes in $\mathcal{M}$. Accordingly, those edges $\epsilon_1$ and $\epsilon_2$ which produce the exact same pattern of jumps and at the exact same order of derivatives in the output responses of the nodes in $\mathcal{M}$ would satisfy $\mathcal{I}$ $(\mathcal{M},\epsilon_1)$ $=$ $\mathcal{I}$ $(\mathcal{M},\epsilon_2)$; and none of them can be identified using just the nodes in $\mathcal{M}$. We further define two set functions $f_D$ and $f_I$ which take a subset of nodes $\mathcal{M}$ and map it to the number of edges that are, respectively, not detectable or not isolatable using the sensor set $\mathcal{M}$. In \cite{rahimian2014detection} these functions are shown to be supermodular; wherefore per the theory of submodular set coverings \cite{Wolsey1982}, adding nodes greedily with respect to these functions would guarantee that the chosen sensor set is within a factor $\log(|\mathcal{E}|)$ of the minimal sensor sets that achieve  detection or isolation goals (solutions of Problems \ref{prob:detection} and \ref{prob:isolation}). The following algorithms are proposed in \cite{CDC13}, and included below for completeness, to implement this idea of supermodular greedy minimization.
\begin{algorithm}
\caption{Determine a Solution $\mathcal{M}_D$ to Problem \ref{prob:detection}}
\label{routine:detection}
\begin{algorithmic}[1]
\REQUIRE $\mathcal{G} = (\mathcal{V},\mathcal{E})$
\State $\mathcal{M}_D \Leftarrow \varnothing$
\WHILE{$f_D(\mathcal{M}_D) \neq 0$}
\STATE $\nu_q  \Leftarrow \arg\min \{ f_D({\mathcal{M}}_D\cup\{\nu_q\}) - f_D({\mathcal{M}}_D);\nu_q \in \mathcal{V}\fgebackslash {\mathcal{M}}_D\}$
\STATE ${\mathcal{M}}_D \Leftarrow {\mathcal{M}}_D\cup\{\nu_q\}$
\ENDWHILE
\ENSURE ${\mathcal{M}}_D$
\end{algorithmic}
\end{algorithm} 

\begin{algorithm}
\caption{Determine a Solution $\mathcal{M}_I$ to Problem \ref{prob:isolation}}
\label{routine:isolation}
\begin{algorithmic}[1]
\REQUIRE $\mathcal{G} = (\mathcal{V},\mathcal{E}) \And \mathcal{M}_D$
\State $\mathcal{M}_I \Leftarrow \mathcal{M}_D$
\WHILE {$f_I(\mathcal{M}_I) \neq 0 \And \mathcal{M}_I \neq \mathcal{V}$}
\STATE $\nu_q  \Leftarrow \arg\min \{ f_I({\mathcal{M}}_I\cup\{\nu_q\}) - f_I({\mathcal{M}}_I);\nu_q \in \mathcal{V}\fgebackslash {\mathcal{M}}_I\}$
\STATE ${\mathcal{M}}_I \Leftarrow {\mathcal{M}}_I\cup\{\nu_q\}$
\ENDWHILE
\IF{$f_I(\mathcal{M}_I) \neq 0$} \State $\mathcal{M}_I \Leftarrow \varnothing$
\ENDIF 
\ENSURE ${\mathcal{M}}_I$
\end{algorithmic}
\end{algorithm} 

It was noted in  Subsection~\ref{sec:bidirectional} that the set of bidirectional links $\mathcal{B}\subset\mathcal{E}$ should be made known to the designer. To facilitate the application of Algorithms~\ref{routine:detection}~and~\ref{routine:isolation} to the case of bidirectional link failures in networked single-integrator agents, we define an equivalence relation $\sim$ on the set $\mathcal{E}$ that identifies two parallel edges in reverse directions only if they are bidirectional. Specifically, for any $\{\tau,\nu\}\subset\mathcal{V}$ such that $(\tau,\nu)\in\mathcal{B}$ and $(\nu,\tau)\in\mathcal{B}$, set $(\tau,\nu)\sim(\nu,\tau)$, while for any two edges $\{\alpha,\beta\}\subset\mathcal{E}$ that $\{\alpha,\beta\}\cap\mathcal{E}\fgebackslash\mathcal{B}\neq\varnothing$, $\alpha\sim\beta$ iff $\alpha=\beta$. The task of the equivalence relation $\sim$ is to identify those edges who share the same head and tail
but at opposite directions, only if they are bidirectional. Every other edge in the network is distinguished and therefore identified only with itself. With the afore-defined equivalence relation $\sim$,
for any edge $\epsilon\in\mathcal{E}$, $\llbracket\epsilon\rrbracket=\{\hat{\epsilon}\in\mathcal{E};\hat{\epsilon}\sim\epsilon\}$ denotes the equivalence class of $\epsilon$, and for any subset of
edges $\mathcal{X}\subset\mathcal{E}$, $\left(\bigslant{\mathcal{X}}{\sim}\right)$
is the quotient of $\mathcal{X}$ by $\sim$, which is the set of all equivalence classes of the elements of $\mathcal{X}$. Last but not least, is the issue of self-loops which are specific to the case of single-integrator agents. In particular, every self-loop would always satisfy an $\mathcal{R}_0$ relation with all nodes in the network and the proposed algorithms cannot be applied for the detection and isolation of a self-loop $(\nu_i,\nu_i)$, although its value (weight) is allowed to change with the failure of a link incoming to node $\nu_i$. In the sequel, the set of all self-loops in $\mathcal{G}$ is denoted by $\mathcal{H}$. Next, changing the definitions of the correspondence $\mathcal{I}(\cdot)$, and the supermodular functions $f_D(\cdot)$ and $f_I(\cdot)$ as follows, allows us to apply Algorithms~\ref{routine:detection}~and~\ref{routine:isolation} to single-integrator networks, while properly identifying bidirectional links and accommodating self-loops. Define for all $\mathcal{M} \subset \mathcal{V}$ and any of the equivalence classes in $\bigslant{\mathcal{E}\fgebackslash\mathcal{H}}{\sim}$:
\begin{align}
&f_D: \mathcal{P}(\mathcal{V}) \to \mathbb{N}_{|\mathcal{E}\fgebackslash\mathcal{H}|}\cup\{0\}, \\ 
& f_D(\mathcal{M}) = |\{\iota \in \left(\bigslant{\mathcal{E}\fgebackslash\mathcal{H}}{\sim}\right) : \forall p \in {\mathcal{M}}, \forall \epsilon \in \iota, (p,\epsilon) \in \mathcal{R}_{0}\}|, \nonumber\\
&\mathcal{I}: \mathcal{P}(\mathcal{V})\times \left(\bigslant{\mathcal{E}\fgebackslash\mathcal{H}}{\sim}\right) \to \mathcal{P}({(\mathbb{N}_z\cup\{0\}) \times \mathcal{V}}), \\ 
&\mathcal{I}(\mathcal{M},\llbracket {\epsilon}\rrbracket ) = \{ (k,p) \in (\mathbb{N}_z\cup\{0\}) \times \mathcal{M}: (p,\epsilon) \in \mathcal{R}_{k} \},\nonumber \\
&f_I: \mathcal{P}(\mathcal{V}) \to \mathbb{N}_{|\mathcal{E}\fgebackslash\mathcal{H}|}\cup\{0\}, \;  f_I(\mathcal{M}) =  \\
& |\{\iota \in  \left(\bigslant{\mathcal{E}\fgebackslash\mathcal{H}}{\sim}\right)  : \exists \hat{\epsilon} \in\mathcal{E}\fgebackslash(\mathcal{H}), \hat{\epsilon} \not\in \iota, \mathcal{I}(\mathcal{M},\iota) = \mathcal{I}(\mathcal{M},\llbracket{\hat{\epsilon}}\rrbracket)\}|.  \nonumber 
\end{align}

\section{Computer Experiments with Large Networks} \label{sec:examples}

In the following subsections, the performance of the developed routines is tested for different random graph models and with varying model parameters.

\subsection{A Random Geometric Graph}\label{geometric}
\begin{figure}[ht]
\centering
\subfigure[Bidirectional Links]{
\includegraphics[width = 0.45\textwidth]{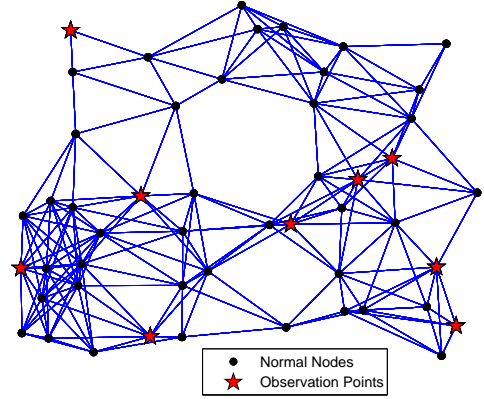}
\label{fig:UndirectedandRandomGeometric}
}
\subfigure[Pairs of Unidirectional Links in Opposite Directions]{
\includegraphics[width = 0.45\textwidth]{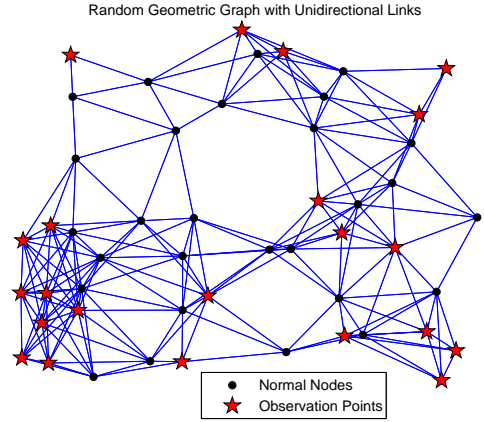}
\label{fig:UnidirectionalRandomGeometric}
}
\caption{\subref{fig:UndirectedandRandomGeometric} With nine observation points and by observing order of derivatives upto $z = 7$, every link in this network is detectable, but none of them can be isolated. \subref{fig:UnidirectionalRandomGeometric} When the undirected links are regarded as pairs of unidirectional links, the required number of observation points increases to $22$ nodes.}
\end{figure}

In a random geometric graph model the nodes of the network are randomly  and uniformly spread across a bounded region, and there is an undirected edge between a pair of nodes, wherever a certain distance threshold is met. The graph of Fig.~\ref{fig:UndirectedandRandomGeometric} depicts one such graph instance with $50$ nodes and $200$ undirected edges, which are interpreted as pairs of bidirectional links. For this graph a total of nine nodes is sufficient for complete detection, whereas even with all of the nodes observed none of the bidirectional links can be isolated. In other words, for any bidirectional link in the network, there exists at least one other link whose removal will induce the same set of jumps in the entire node-set of the network.

The situation is rather different if the $200$ undirected edges of the network in Fig.~\ref{fig:UndirectedandRandomGeometric} are regarded as $400$ unidirectional links. Then the output of Routine~\ref{routine:detection} has $22$ nodes that are indicated in Fig.~\ref{fig:UnidirectionalRandomGeometric}, and by observing them the designer can isolate $280$ edges out of the total $400$. Observing all of the nodes in the network decreases the cardinality of the set of unresolved edges from $120$ to just $93$, out of the total $400$. It is worth highlighting that with the change in the interpretation of the links from bidirectional to unidirectional, matrix $R$ of the network remains the same, and so does the required highest order of derivatives $z = 7$.

\begin{figure}[ht]
\centering
\includegraphics[width = 0.48\textwidth]{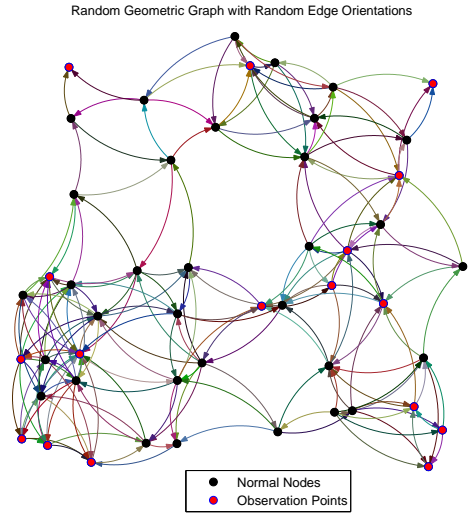}
\caption{If the edges of the undirected random geometric graph are randomly oriented, then $17$ nodes would be enough to achieve the detection task. This also increases the required highest order of derivatives from $z = 7$ to $z = 9$.}
\label{fig:OrientedRandomGeometric}
\end{figure}

\begin{figure}[ht]
\centering
\includegraphics[width = 0.45\textwidth]{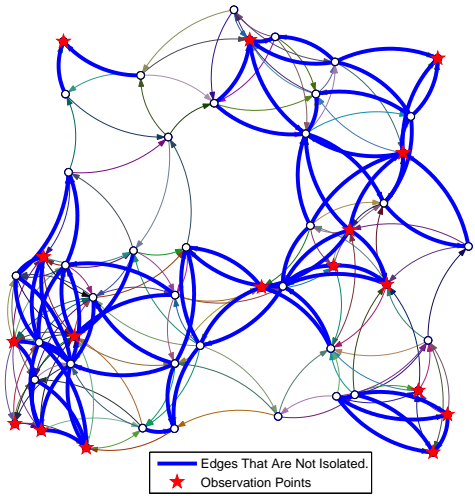}
\caption{The $75$ edges that are highlighted cannot be isolated using the indicated observation points. However, even with all of the nodes observed, there still remain $34$ edges that cannot be isolated.}
\label{fig:OrientedRandomGeometricNOTisolated}
\end{figure}

Next, each of the undirected edges in Fig.~\ref{fig:UndirectedandRandomGeometric} is oriented randomly leading to a total of $200$ unidirectional edges in Fig.~\ref{fig:OrientedRandomGeometric}. In the latter, a total of $17$ nodes is sufficient for detection, and these $17$ nodes enable the isolation of all but $75$ edges of the digraph, which are highlighted in Fig.~\ref{fig:OrientedRandomGeometricNOTisolated}. For this directed network, by observing all of the nodes in the network, the cardinality of the set of unresolved links reduces to $34$.

The preceding results suggest that while detection is achievable more easily in undirected networks, the increased diversity brought about by the directionality of the links improves the isolation task for the case of directed networks.

The focus of investigation in the following subsections is shifted to the Erd\H{o}s-R\'{e}nyi random graph model, for which the role of edge probability and graph size on the cardinality of the detection set and the highest order of derivatives required is explored.

\subsection{Erd\H{o}s-R\'{e}nyi Random Graphs: Directed versus Undirected}\label{erdos}

\begin{figure}
\centering
\includegraphics[width = 0.39\textwidth]{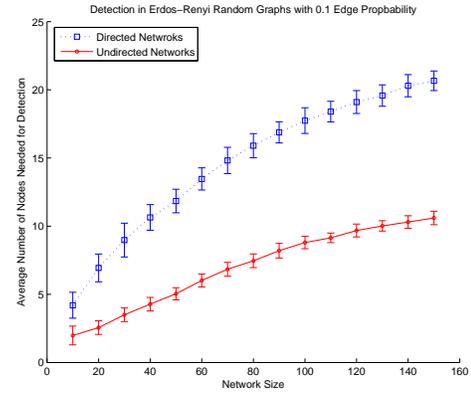}
\caption{$|\mathcal{M}_D|$ versus $n$ for  Erd\H{o}s-R\'{e}nyi random graphs with $p = 0.1$}
\label{fig:erdosNumberVersusSize}
\end{figure}

\begin{figure}
\centering
\includegraphics[width = 0.39\textwidth]{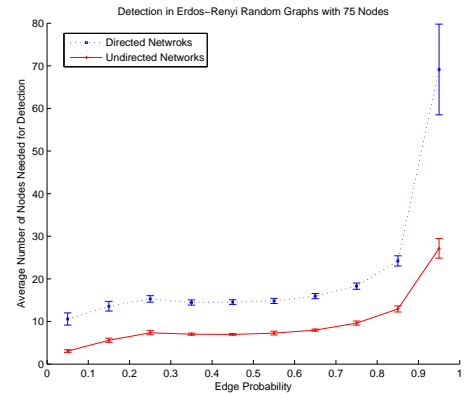}
\caption{$|\mathcal{M}_D|$ versus $p$ for  Erd\H{o}s-R\'{e}nyi random graphs with $n = 75$}
\label{fig:erdosNumberVersusProbability}
\end{figure}

In a Erd\H{o}s-R\'{e}nyi random graph model every potential edge is either existent or not with a fixed probability $p$, and independently of all the rest. This model is implemented for varying network sizes $n$, and different edge probabilities $p$. In Figs.\ref{fig:erdosNumberVersusSize} and \ref{fig:erdosNumberVersusProbability}, the cardinality of the detection sets $|\mathcal{M}_D|$ in several randomly generated instances are recorded, averaged, and plotted. The sample means in each case are computed over $50$ random instances and the error bars indicate the sample standard deviations for those instances. The plots in all cases confirm the increased difficulty of the detection process for the case of directed networks. Moreover, the cardinality of the detection sets does not scale fast with the network size; an observation which is of practical significance for large networks and complements the theoretical guarantees that are available from the submodular set covering literature. In the case of edge probabilities, however, it is observed that as the edge probabilities approach $1$ leading to a complete graph, the number of nodes required for detection becomes increasingly large. It is worth highlighting that although for small $p$ the networks are sparse and can have large diameters, as the edge probability is increased beyond $0.3$ the network diameter remains constant at $2$ so that only the first three derivatives of any chosen sensor set need to be observed. Similarly for $n$, as the network size is increased beyond $60$, the network diameters remain fixed at $3$ and only the first four derivatives of the outputs in any chosen sensor set are observed.

\section{Conclusions}\label{sec:conc}

In this paper, we developed FDI techniques for single-integrator networks that enable the designer to detect and isolate link failures based on the observed jumps in the derivatives of the output responses of a subset of nodes by relating the jumps in the derivatives at the time of failure to the distance of the failed link from the observation point. Our results covered both cases of unidirectional and bidirectional link failures. We also extended our previously developed sensor placement algorithms to accommodate both types of link failures (unidirectional and bidirectional). These algorithms were tested in large random networks, and the results suggest that link failures in directed networks are harder to detect but easier to isolate, as compared to undirected networks. The latter effect can be attributed to the increased diversity that is brought about by the directionality of the links. Moreover, both the cardinality of the detection sets and the required order of derivatives are shown to scale up reasonably well with the network size, and this agrees with the performance guarantees that are available from the theory of submodular set coverings and bound the size of the chosen sets to within a multiplicative $\log(|\mathcal{E}|)$ factor of the minimal sensor set, where $|\mathcal{E}|$ is the size of the edge set.

\appendix[Proofs of The Main Results]

\subsection{Theorem~\ref{theo:detection2}}

Given an initial condition $\mathbf{x}_{0}:=\mathbf{x}(t_{0})\in\mathbb{R}^{N}$,
the solution to \eqref{eq:laplaciandynamics} for $t_{0}\leq t\leq t_{f}$
is trivially given by: 
\begin{equation}
\mathbf{x}(t)=e^{A\left(t-t_{0}\right)}\mathbf{x}_{0}+\int_{t_{0}}^{t}\mathrm{e}^{A(t-\tau)}B\mathbf{w}(\tau)\,\mathrm{d}\tau,\mbox{ for }t_{f}\geqslant t\geqslant t_{0}.\label{eq:matrixsolutionlaplaciandynamics}
\end{equation}
The evolution of states after failure is therefore governed by the
state matrix $\bar{A}$, instead of $A$, as follows
\begin{equation}
\mathbf{x}\left(t\right)=\left(e^{\bar{A}(t-t_{f})}\mathbf{x}_{f}+\int_{t_{f}}^{t}\mathrm{e}^{\bar{A}(t-\tau)}B\mathbf{w}(\tau)\,\mathrm{d}\tau\right),\mbox{ for }t\geqslant t_{f},\label{eq:agentlaplaciandynamicsFAILURE}
\end{equation}
where $\mathbf{x}_{f}:=\mathbf{x}(t_{f}^{-})$, i.e., the state of
the faultless evolution at the instant right before failure.

For any $p\in\mathbb{N}_{N}$ fixed, differentiating \eqref{eq:matrixsolutionlaplaciandynamics}
and \eqref{eq:agentlaplaciandynamicsFAILURE} $k$ times and using
the Leibniz integral rule yields for $t_{f}>t>t_{0}$ that: 
\begin{align}
&\frac{{d}^{k}}{{d}t^{k}}x_{p}(t)  = \label{eq:derivates1} \\ &\mathbf{e}_{p}^{T}\left(A^{k}\mathbf{x}(t)+\sum_{m=0}^{k-2}A^{m}B\frac{{d}^{k-m-1}}{{d}t^{k-m-1}}\mathbf{w}(t)+A^{k-1}B\mathbf{w}(t)\right),\nonumber
\end{align}
and for $t>t_{f}$ that,
\begin{align}
&\frac{{d}^{k}}{{d}t^{k}}x_{p}(t)  = \label{eq:derivates2} \\ 
& \mathbf{e}_{p}^{T}\left(\bar{A}^{k}\mathbf{x}(t)+\sum_{m=0}^{k-2}\bar{A}^{m}B\frac{{d}^{k-m-1}}{{d}t^{k-m-1}}\mathbf{w}(t)+\bar{A}^{k-1}B\mathbf{w}(t)\right).\nonumber
\end{align}
Next, note that by differentiability of $\mathbf{w}(t)$ and continuity
of the states, $\frac{{d}^{k}}{{d}t^{k}}\mathbf{w}(t)\big|_{t=t_{f}^{+}}$
$=$ $\frac{{d}^{k}}{{d}t^{k}}\mathbf{w}(t)_{t=t_{f}^{-}}$
$=$ $\frac{{d}^{k}}{{d}t^{k}}\mathbf{w}(t)\big|_{t=t_{f}},\forall k\in\mathbb{N}$
and $\mathbf{x}(t_{f}^{+})=\mathbf{x}(t_{f}^{-})$. Hence, subtracting
the two equations in \eqref{eq:derivates1} and \eqref{eq:derivates2} for $t_{f}^{+}$ and $t_{f}^{-}$
yields:

\begin{align}
&\Delta_{t_f}({p,k})=  \mathbf{e}_{p}^{T}(\bar{A}^{k}-A^{k})\mathbf{x}(t_{f})+ \\&  \mathbf{e}_{p}^{T}(\bar{A}^{k-1}-A^{k-1})B\mathbf{w}(t_{f})+ \\& \mathbf{e}_{p}^{T}\sum_{m=0}^{k-2}(\bar{A}^{m}-A^{m})B\left(\frac{{d}^{k-m-1}}{{d}t^{k-m-1}}\mathbf{w}(t)\bigg|_{t=t_{f}}\right),\label{mainEQ}
\end{align}

With $p$, $\bar{A}$ and $A$ fixed in the preceding and $k\in\mathbb{N}_{\textrm{dist}(\nu_{j},\nu_{p})}$,
for all $m\in\mathbb{N}_{\textrm{dist}(\nu_{j},\nu_{p})}$,
define $\nabla(q,m)=\left[\bar{A}^{m}\right]_{pq}-\left[A^{m}\right]_{pq}$,
so that \eqref{mainEQ} can be rewritten as: 
\begin{align}
&\Delta_{t_f}({p,k})=\sum_{q=1}^{N}\nabla(q,k){x}_{q}(t_{f})+\sum_{q=1}^{N}\nabla(q,k-1)\left[B\mathbf{w}(t_{f})\right]_{q} \nonumber \\ & + \sum_{m=0}^{k-2}\sum_{q=1}^{N}\nabla(q,m)\left[B\frac{{d}^{k-m-1}}{{d}t^{k-m-1}}\mathbf{w}(t_{f})\right]_{q}\label{mainEQ2}
\end{align}
Next to compute $\nabla(q,m)$, substitute for $\left[A^{m}\right]_{pq}$
and $\left[\bar{A}^{m}\right]_{pq}$ from \eqref{eq:numWalks} to
get: 
\begin{equation}
\nabla(q,m)=  \label{eq:sumWeights1}\\ \Phi(\Omega^{m}(\bar{\mathcal{G}};{\nu_{q},\nu_{p}}),\bar{A})-\Phi(\Omega^{m}({\mathcal{G}};{\nu_{q},\nu_{p}}),{A}). \nonumber
\end{equation}
By partitioning the sets $\Omega^{m}(\bar{\mathcal{G}};{\nu_{q},\nu_{p}})$
and $\Omega^{m}({\mathcal{G}};{\nu_{q},\nu_{p}})$, \eqref{eq:sumWeights1}
can be rewritten as: 
\begin{align}
\nabla(q,m)= & [\Phi(\Omega^{m}(\bar{\mathcal{G}};{\nu_{q},\nu_{p}})\fgebackslash\Omega^{m}(\bar{\mathcal{G}};{\nu_{q},\nu_{i},\nu_{p}}),\bar{A}) \nonumber \\ & +\Phi(\Omega^{m}(\bar{\mathcal{G}};{\nu_{q},\nu_{i},\nu_{p}}),\bar{A})]-\nonumber \\
 & [\Phi(\Omega^{m}(\mathcal{G};\nu_{q},\nu_{p})\fgebackslash\Omega^{m}(\mathcal{G};{\nu_{q},\nu_{i},\nu_{p}}),A) \nonumber \\ &  + \Phi(\Omega^{m}(\mathcal{G};{\nu_{q},\nu_{i},\nu_{p}}),A)].\label{eq:partitioned}
\end{align}

Next note that none of the walks in $\Omega^{m}(\mathcal{G};\nu_{q},\nu_{p})$ $\fgebackslash$ $\Omega^{m}(\mathcal{G};{\nu_{q},\nu_{i},\nu_{p}})$ or $\Omega^{m}(\bar{\mathcal{G}};{\nu_{q},\nu_{p}})$ $\fgebackslash$ $\Omega^{m}(\bar{\mathcal{G}};{\nu_{q},\nu_{i},\nu_{p}})$ can include $(\nu_{\gamma},\nu_{i})$ as an edge for any $\gamma\in\mathbb{N}_{N}$. This is true, since otherwise if there exists a walk $\mathcal{W}_{1}$
that violates the above, then removing the segment of $\mathcal{W}_{1}$
from $\nu_{q}$ to $\nu_{i}$, which consists of at least two edges,
one to reach $\nu_{\gamma}$ from $\nu_{q}$ followed by the edge
$(\nu_{\gamma},\nu_{i})$, yields a new $\nu_{i}\nu_{p}$ walk $\mathcal{W}_{2}$
with length at most $m-2$. Now $(\nu_{j},\nu_{i})\mathcal{W}_{2}$
is a $\nu_{j}\nu_{p}$ walk of length at most $m-1$ in $\mathcal{G}$,
which is a contradiction, since $m\leqslant\textrm{dist}(\nu_{j},\nu_{p})$.
It next follows that $\Phi(\Omega^{m}(\mathcal{G};\nu_{q},\nu_{p})$ $\fgebackslash$
$\Omega^{m}(\mathcal{G};{\nu_{q},\nu_{i},\nu_{p}}),$ $A)$ $=$ $\Phi(\Omega^{m}(\bar{\mathcal{G}};{\nu_{q},\nu_{p}})$
$\fgebackslash$ $\Omega^{m}(\bar{\mathcal{G}};{\nu_{q},\nu_{i},\nu_{p}}),$
$\bar{A})$, as none of the walks involved include any of the edges
$(\nu_{q},\nu_{i})$ for $q\in\mathbb{N}_{N}$, and these are the
only edges at which the digraphs $\bar{\mathcal{G}}$ and $\mathcal{G}$
or the in-weightings $\bar{A}$ and $A$ differ. Hence, \eqref{eq:sumWeights1}
simplifies into:

\begin{align}
&\nabla(q,m)  = \label{eq:simplified1} \\ &\Phi(\Omega^{m}(\bar{\mathcal{G}};{\nu_{q},\nu_{i},\nu_{p}}),\bar{A})-\Phi(\Omega^{m}(\mathcal{G};{\nu_{q},\nu_{i},\nu_{p}}),A)\nonumber \\
&  =\bar{a}_{iq}\Phi(\Omega^{m-1}(\bar{\mathcal{G}};\nu_{i},\nu_{p}),\bar{A})-a_{iq}\Phi(\Omega^{m-1}(\mathcal{G};\nu_{i},\nu_{p}),A). 
\end{align}

The last step in deriving a simplified expression for $\nabla(q,m)$
is to argue that $\Phi(\Omega^{m-1}(\mathcal{G};{\nu_{i},\nu_{p}}),$
$A)=\Phi(\Omega^{m-1}(\bar{\mathcal{G}};{\nu_{i},\nu_{p}}),$ $\bar{A})$.
To see why, note that since $\bar{\mathcal{G}}$ is derived upon removal
of the edge $(\nu_{j},\nu_{i})$ from $\mathcal{G}$, it follows that
\begin{align}
\mbox{{dist}}({\nu}_{i},{\nu_{p}}) = \mbox{{dist}}(\bar{\mathcal{G}};{\nu}_{i},{\nu_{p}})
\geqslant \textrm{dist}(\nu_{j},\nu_{p})-1, 
\label{eq:paths}
\end{align} where the digraph argument for the distance function indicate that the distances are calculated with respect to the edge-removed digraph $\bar{\mathcal{G}}$ as opposed to usual case where the distances are calculated with respect to the original digraph $\mathcal{G}$. The inequalities in \eqref{eq:paths} together with $\textrm{dist}(\nu_{j},\nu_{p})\geqslant m$,
implies that $\mbox{{dist}}({\nu}_{i},{\nu_{p}})$ $=$ $\mbox{{dist}}(\bar{\mathcal{G}};{\nu}_{i},{\nu_{p}})$
$\geqslant$ $m-1$, so that none of the walks in $\Omega^{m-1}(\mathcal{G};{\nu_{i},\nu_{p}})$
or $\Omega^{m-1}(\bar{\mathcal{G}};{\nu_{i},\nu_{p}})$ can include
any $(\nu_{q},\nu_{i})$ edges, $\forall q\in\mathbb{N}_{N}$, and
$\Phi(\Omega^{m-1}(\mathcal{G};{\nu_{i},\nu_{p}}),$ $A)=\Phi(\Omega^{m-1}(\bar{\mathcal{G}};{\nu_{i},\nu_{p}}),$
$\bar{A})$; hence, 
\begin{align}
\nabla(q,m) & =\left(\bar{a}_{iq}-a_{iq}\right)\Phi(\Omega^{m-1}(\mathcal{G};\nu_{i},\nu_{p}),A),\label{eq:simplified2}
\end{align}
which upon replacement in \eqref{mainEQ2} yields: 
\begin{align}
&{\Delta}_{t_f}(p,k)=  \Phi(\Omega^{k-1}(\mathcal{G};\nu_{i},\nu_{p}),A)\sum_{q=1}^{N}\left(\bar{a}_{iq}-a_{iq}\right){x}_{q}(t_{f})+ \nonumber \\& \Phi(\Omega^{k-2}(\mathcal{G};\nu_{i},\nu_{p}),A)\sum_{q=1}^{N}\left(\bar{a}_{iq} - a_{iq}\right)\left[B\mathbf{w}(t_{f})\right]_{q} \nonumber \\
 & +\sum_{m=1}^{k-2}\Phi(\Omega^{m-1}(\mathcal{G};\nu_{i},\nu_{p}),A)\sum_{q=1}^{N}\left(\bar{a}_{iq}-a_{iq}\right) \nonumber \\& \left[B\frac{{d}^{k-m-1}}{{d}t^{k-m-1}}\mathbf{w}(t_{f})\right]_{q}.\label{mainfinal1}
\end{align}
To complete the proof, note that for $m\leqslant\textrm{dist}(\nu_{j},\nu_{p})-1$,
$\Omega^{m-1}(\mathcal{G};\nu_{i},\nu_{p})=\varnothing$, since for
any $\mathcal{W}_{3}\in\Omega^{m-1}(\mathcal{G};\nu_{i},\nu_{p})$,
$(\nu_{j},\nu_{i})\mathcal{W}_{3}$ is a $\nu_{j}\nu_{p}$ walk of
length $m$ which contradicts with $m\leqslant\textrm{dist}(\nu_{j},\nu_{i})-1$.
Thence, for $k\leqslant\textrm{dist}(\nu_{j},\nu_{p})$,
\eqref{mainfinal1} simplifies into ${\Delta}_{t_f}(p,k)=\Phi(\Omega^{k-1}(\mathcal{G};\nu_{i},\nu_{p}),A)\sum_{q=1}^{N}\left(\bar{a}_{iq}-a_{iq}\right){x}_{q}(t_{f})$,
thus completing the proof for the case of the failure of the single
link $\bar{\epsilon}$. That for $k<\textrm{dist}(\nu_{j},\nu_{p})$,
${\Delta}_{t_f}(p,k)=0$ also follows as $\Omega^{k-1}(\mathcal{G};\nu_{i},\nu_{p})=\varnothing$
for any such $k$. \hfill{}{\scriptsize $\blacksquare$}{\scriptsize \par}

\subsection{Lemma~\ref{lem:singleIntegratorDetection}}

Notice that any shortest path from $\nu_{i}$ to $\nu_{p}$ of length $k-1$ gives a path of length $k$ from $\nu_j$ to $\nu_p$, whence given $k=\mbox{{dist}}(\nu_{j},\nu_{p})$, it follows that:
\begin{align}
\mbox{{dist}}(\nu_{i},\nu_{p}) \geqslant k-1.
\label{cond:GEQ} 
\end{align}
Also notice that if $\mbox{{dist}}(\nu_{i},\nu_{p})\leqslant$ $k-1$ there is no path of length ${k-1}$ from $\nu_{i}$ to $\nu_{p}$, i.e. 
\begin{align}
& \Phi(\Omega^{k-1}(\mathcal{G};\nu_{i},\nu_{p}),A)\mbox{ is nonzero,} \\ &\mbox{only if } \mbox{{dist}}(\nu_{i},\nu_{p})\leqslant k-1.
\label{cond:PHI}
\end{align}
Now \eqref{cond:GEQ} and \eqref{cond:PHI} together imply that for the right-hand side of \eqref{eq:mainResult} in Theorem~\ref{theo:detection2} to be non-zero when $k=\mbox{{dist}}(\nu_{j},\nu_{p})$, it should be true that $\mbox{{dist}}(\nu_{i},\nu_{p})+1$ $=$ $k$, which is the same as the claimed condition.   \hfill{}{\scriptsize $\blacksquare$}{\scriptsize \par}

\subsection{Proposition~\ref{prop:detection3}}

To see how Theorem~\ref{theo:detection2} applies to the case of bidirectional link failures, in parallelism with $\bar{A}$ and $\bar{\mathcal{G}}$ in the preceding, let $\hat{A} \neq A$ be an in-weighting on $\hat{\mathcal{G}}$ that denotes the perturbed version of $A$ following the simultaneous failure of links $\bar{\epsilon}$ and $\hat{\epsilon}$. The perturbations only affect the entries of $A$ on its $i-$th and $j-$th rows, such that $\hat{a}_{ij} = \hat{a}_{ji} = 0$, while $\hat{a}_{qr} = {a}_{qr}$, $\forall r \in \mathbb{N}_N$ and $\forall q \in \mathbb{N}_N\fgebackslash \{i,j\}$. For any agent $p \in \mathbb{N}_N$ and $t \geqslant t_f$ the evolution of the state of agent $p$ following the simultaneous failure of links $\bar{\epsilon}$ and $\hat{\epsilon}$ is given by \eqref{eq:agentlaplaciandynamicsFAILURE} with $\hat{A}$ substituted for $\bar{A}$. Repeating the same procedure as in the proof of Theorem~\ref{theo:detection2} leads to the proof of the proposition as follows.

For the case of simultaneous failure of $\bar{\epsilon}$ and $\hat{\epsilon}$,
note that \eqref{eq:derivates1} to \eqref{eq:partitioned} continue
to hold after substituting $\hat{A}$ and $\hat{\mathcal{G}}$ for
$\bar{A}$ and $\bar{\mathcal{G}}$, respectively. The transitions
from \eqref{eq:partitioned} to \eqref{eq:simplified1} and \eqref{eq:simplified2}
also carry through with the same replacements and upon the additional
observation that the walks in $\Omega^{m}(\mathcal{G};\nu_{q},\nu_{p})\fgebackslash\Omega^{m}(\mathcal{G};{\nu_{q},\nu_{i},\nu_{p}})$,
$\Omega^{m}(\hat{\mathcal{G}};{\nu_{q},\nu_{p}})\fgebackslash\Omega^{m}(\hat{\mathcal{G}};{\nu_{q},\nu_{i},\nu_{p}})$,
$\Omega^{m-1}(\mathcal{G};{\nu_{i},\nu_{p}})$ and $\Omega^{m-1}(\hat{\mathcal{G}};{\nu_{i},\nu_{p}})$
include neither any $(\nu_{q},\nu_{i})$ edges, as stated in the previous
case, nor any $(\nu_{q},\nu_{j})$ edges, $\forall q\in\mathbb{N}_{N}$.
The rest of the proof is identical to the previous case, except for
$\bar{A}$ which should be replaced with $\hat{A}$. \hfill{}{\scriptsize $\blacksquare$}{\scriptsize \par}

\bibliographystyle{IEEEtran}
\bibliography{refDistinguishability,newRef}

\begin{thebibliography}{10}
\providecommand{\url}[1]{#1}
\csname url@samestyle\endcsname
\providecommand{\newblock}{\relax}
\providecommand{\bibinfo}[2]{#2}
\providecommand{\BIBentrySTDinterwordspacing}{\spaceskip=0pt\relax}
\providecommand{\BIBentryALTinterwordstretchfactor}{4}
\providecommand{\BIBentryALTinterwordspacing}{\spaceskip=\fontdimen2\font plus
\BIBentryALTinterwordstretchfactor\fontdimen3\font minus
  \fontdimen4\font\relax}
\providecommand{\BIBforeignlanguage}[2]{{%
\expandafter\ifx\csname l@#1\endcsname\relax
\typeout{** WARNING: IEEEtran.bst: No hyphenation pattern has been}%
\typeout{** loaded for the language `#1'. Using the pattern for}%
\typeout{** the default language instead.}%
\else
\language=\csname l@#1\endcsname
\fi
#2}}
\providecommand{\BIBdecl}{\relax}
\BIBdecl

\bibitem{mesbahiBook}
M.~Mesbahi and M.~Egerstedt, \emph{{Graph Theoretic Methods in Multiagent
  Networks}}.\hskip 1em plus 0.5em minus 0.4em\relax Princeton University
  Press, 2010.

\bibitem{aminAutomatica}
M.~A. Rahimian and A.~G. Aghdam, ``Structural controllability of multi-agent
  networks: Robustness against simultaneous failures,'' \emph{Automatica},
  vol.~49, no.~11, pp. 3149 -- 3157, 2013.

\bibitem{Kleinberg}
J.~Kleinberg, M.~Sandler, and A.~Slivkins, ``Network failure detection and
  graph connectivity,'' \emph{SIAM Journal on Computing}, vol.~38, no.~4, pp.
  1330--1346, 2008.

\bibitem{largeScaleNonlinearPowerNetworks}
W.~Pan, Y.~Yuan, H.~Sandberg, J.~Gon\c{c}alves, and G.-B. Stan, ``Real-time
  fault diagnosis for large-scale nonlinear power networks,'' in
  \emph{Proceedings of the 52nd IEEE Conference on Decision and Control}, 2013,
  pp. 2340--2345.

\bibitem{6545301}
F.~Pasqualetti, F.~Dorfler, and F.~Bullo, ``Attack detection and identification
  in cyber-physical systems,'' \emph{Automatic Control, IEEE Transactions on},
  vol.~58, no.~11, pp. 2715--2729, 2013.

\bibitem{discreteTimeFaultDetection}
E.~E. Tiniou, P.~M. Esfahani, and J.~Lygeros, ``Fault detection with
  discrete-time measurements: An application for the cyber security of power
  networks,'' in \emph{Proceedings of the 52nd IEEE Conference on Decision and
  Control}, 2013, pp. 194--199.

\bibitem{6567892}
P.~Menon and C.~Edwards, ``Robust fault estimation using relative information
  in linear multi-agent networks,'' \emph{Automatic Control, IEEE Transactions
  on}, vol.~PP, no.~99, pp. 1--1, 2013.

\bibitem{6482175}
C.~Keliris, M.~Polycarpou, and T.~Parisini, ``A distributed fault detection
  filtering approach for a class of interconnected continuous-time nonlinear
  systems,'' \emph{Automatic Control, IEEE Transactions on}, vol.~58, no.~8,
  pp. 2032--2047, 2013.

\bibitem{asjc}
M.~A. Rahimian, A.~Ajorlou, and A.~G. Aghdam, ``Digraphs with distinguishable
  dynamics under the multi-agent agreement protocol,'' \emph{Asian Journal of
  Control}, 2014, in press.

\bibitem{CDC13}
M.~A. Rahimian and V.~M. Preciado, ``Detection and isolation of link failures
  under the agreement protocol,'' in \emph{Proceedings of the 52nd IEEE
  Conference on Decision and Control}, 2013, pp. 7364--7369.

\bibitem{BiggsGraphTheory}
N.~Biggs, \emph{{Algebraic Graph Theory}}.\hskip 1em plus 0.5em minus
  0.4em\relax Cambridge University Press, 1994.

\bibitem{Preciado:2013:MSA:2502376.2502379}
V.~M. Preciado and A.~Jadbabaie, ``Moment-based spectral analysis of
  large-scale networks using local structural information,'' \emph{IEEE/ACM
  Transactions on Networking}, vol.~21, no.~2, pp. 373--382, Apr. 2013.

\bibitem{rahimian2014detection}
M.~A. Rahimian and V.~M. Preciado, ``Detection and isolation of failures in
  directed networks of lti systems,'' \emph{arXiv preprint arXiv:1408.3164},
  2014.

\bibitem{Wolsey1982}
L.~Wolsey, ``An analysis of the greedy algorithm for the submodular set
  covering problem,'' \emph{Combinatorica}, vol.~2, no.~4, pp. 385 -- 393,
  1982.

\end{thebibliography}

\end{document}